\documentclass[11pt]{amsart}
\usepackage{amssymb}

\input xy
\xyoption {all}
\input xy
\xyoption {all}

\newcommand{\comment}[1]{}
\newtheorem{theorem}{Theorem}
\newtheorem {lemma}{Lemma}
\newtheorem{conjecture}{Conjecture}

\theoremstyle{definition}
 
\theoremstyle {definition}

\begin{document}
\baselineskip=14.5pt
\title {The combinatorics of Lehn's conjecture}

\author{A. Marian}
\address{Department of Mathematics, Northeastern University}
\email {a.marian@neu.edu}
\author{D. Oprea}
\address{Department of Mathematics, University of California, San Diego}
\email {doprea@math.ucsd.edu}
\author{R. Pandharipande}
\address{Department of Mathematics, ETH Z\"urich}
\email {rahul@math.ethz.ch}

\begin{abstract}Let $S$ be a nonsingular projective surface equipped with a line bundle $H$. Lehn's conjecture is a formula for the top Segre class of the tautological bundle associated to $H$ on the Hilbert scheme of points of $S.$ Voisin has recently reduced Lehn's conjecture to the vanishing of certain coefficients of special power series. The first result here is a proof of the vanishings required by Voisin by residue calculations (A. Szenes and M. Vergne have independently found the same proof). Our second result is an elementary solution of the parallel question for the top Segre class on the symmetric power of a nonsingular projective curve $C$ associated to a higher rank vector bundle $V$ on $C$. Finally, we propose a complete conjecture for the top Segre class on the Hilbert scheme of points of $S$ associated to a higher rank vector bundle on $S$ in the $K$-trivial case. 
\end{abstract}
\maketitle 
\vskip.5in
\noindent{\bf Lehn's conjecture.} The number of $(n-2)$-subspaces in $\mathbb P^{2n-2}$ which are $n$-secant to a nonsingular curve
$$C\subset \mathbb P^{2n-2}$$ 
of genus $g$ and degree $d$
is a classical
enumerative calculation \cite {ACGH}. The answer can be expressed in terms of Segre integrals on the symmetric
{\footnote{The $n^{th}$ symmetric product of $C$ is the Hilbert scheme of points $C^{[n]}$. For curves $C$ and surfaces $S$,  we use the standard notation for the tautological bundle $H^{[n]}$ of rank $n$ on the Hilbert schemes $C^{[n]}$ and $S^{[n]}$ associated to a line bundle $H$, see [EGL].}}
 product
$C^{[n]}$ of $C$.
Let the line bundle
$$H\to C$$ be the degree $d$ restriction of $\mathcal{O}_{\mathbb P^{2n-2}}(1)$.
The $n$-secant problem is solved by the  Segre integral, and the answer can be written in closed form \cite {LeB}, \cite {C},
\begin{equation}
\label{curve}\sum_{n=0}^{\infty} z^n\int_{C^{[n]}} s_{n}({H}^{[n]})=\frac{(1-w)^{d+2\chi(\mathcal O_C)}}{(1-2w)^{\chi (\mathcal O_C)}}\, ,\end{equation}
after the change of variables $$z=w(1-w)\, .$$ 

Going further, consider a pair $(S, H)$ consisting of a nonsingular projective surface and a line bundle $H\to S$. 
The Segre integrals
 $$\int_{S^{[n]}} s_{2n}(H^{[n]})$$ 
on the Hilbert scheme of points $S^{[n]}$
count
the $n$-secants of dimension $n-2$ to the image of the surface 
$$S\to \mathbb P^{3n-2}\,, \ \ \ H=\mathcal{O}_{\mathbb P^{3n-2}}(1)|_S\, .$$
 The following conjecture was made by Lehn \cite {L}:
\begin{equation}\label{series}\sum_{n=0}^{\infty} z^n \int_{S^{[n]}} s_{2n}(H^{[n]})=\frac{(1-w)^{a}(1-2w)^{b}}{(1-6w+6w^2)^{c}}\, \end{equation} for constants $$a=H\cdot K_S-2K_S^2\, ,\,\,\ b=(H-K_S)^2+3\chi(\mathcal O_S)\, ,\,
\, \ 
c=\frac{1}{2}H(H-K_S)+\chi(\mathcal O_S)\, .$$ A more complicated change of variables is needed here,
 $$z=\frac{w(1-w)(1-2w)^4}{(1-6w+6w^2)^3}\, .$$ The first few terms are $$z=w+9w^2+68w^3+\ldots\iff w=z-9z^2+94z^3+\ldots\,.$$
\vskip.1in

For $K$-trivial surfaces, Lehn's conjecture was established in \cite {MOP} via a study of the virtual geometry of a suitable Quot scheme. The results in \cite {V} on blowups of $K3$ surfaces, obtained via classical geometry, provide the missing geometric pieces needed to establish Lehn's conjecture in full generality. 

\begin{theorem}\label{tle} Lehn's conjecture holds for all surfaces. 
\end{theorem} 

\noindent {\it Proof.} By the results of \cite {EGL}, the Segre series can be written in the form \begin{equation}\label{split}\sum_{n=0}^{\infty} z^n \int_{S^{[n]}} s_{2n}(H^{[n]})=A_1(z)^{H^2} \cdot A_2(z)^{\chi(\mathcal O_S)} \cdot A_3(z)^{H\cdot K_S} \cdot A_4(z)^{K_S^2}\end{equation} for four universal power series $$A_1, A_2, A_3, A_4\in \mathbb Q[[z]]\, .$$ Lehn's conjecture \eqref{series} consists of
the following evaluations: \begin{eqnarray}\label{lehnlehn}A_1(z)=\frac{1-2w}{(1-6w+6w^2)^{\frac{1}{2}}},\,\,\,  &&A_2(z)=\frac{(1-2w)^3}{1-6w+6w^2}\, ,
\\ \nonumber A_3(z)=\frac{(1-w)(1-6w+6w^2)^{\frac{1}{2}}}{(1-2w)^2},\,\,\, &&A_4(z)=\frac{1-2w}{(1-w)^2}\, 
.
\end{eqnarray}

 The expressions for $A_1, A_2$ in \eqref{lehnlehn} were proven correct in \cite {MOP}. Key to the argument was the closed form evaluation of all Segre integrals over Hilbert schemes of points on $K3$ surfaces, \begin{equation}\label{clo}\int_{S^{[n]}}s_{2n}(H^{[n]})=2^n \binom{\frac{H^2}{2}+2-2n}{n}.\end{equation} 

We show that the results in \cite {V} on blowups of $K3$s give the remaining series $A_3$ and $A_4$. To this end, let $S$ be the blowup of a generic primitively polarized $K3$ surface $(X, L)$ at one point. Define the line bundle $$H=L \otimes E^{-k}$$ on $S$ where $E$ is
 the exceptional line bundle on the blowup. We have
 $$ H \cdot K_S=k\, .$$ The crucial input is provided by Theorem $3$ in \cite {V}, which, in our notation, states\footnote{It would be interesting to see if these Segre vanishings can be obtained also by the methods of \cite {MOP}.}\begin{equation}\label{vv}s_{2n}(H^{[n]})=0 \, \, \, \text{whenever} \, \, \, \chi(H) = 3n-1, \, \,  k=n-1, \text{  or } k=n\, .\end{equation} Proposition 19 in \cite {V} furthermore shows that the vanishings \eqref{vv} uniquely determine the series $A_3, A_4$. The series are determined inductively, coefficient by coefficient. However, the closed form expressions for $A_3, A_4$ stated in \eqref{lehnlehn} were left open in \cite{V}. To complete the proof of Lehn's conjecture, it suffices to show  
$$\text{Coeff}_{z^n} \left [ A_1(z)^{H^2} \cdot A_2(z)^{\chi(\mathcal O_S)} \cdot A_3(z)^{H\cdot K_S} \cdot A_4(z)^{K_S^2} \right ] = 0$$
for $\chi(H) = 3n-1, \, \,  k=n-1, \text{  or } k=n,$ where the series $A_1, A_2, A_3, A_4$ are given by \eqref{lehnlehn}.

We will prove more strongly that whenever $\chi (H) = 3n-1,$
\begin{equation}\label{another}\text{Coeff}_{z^n} \left [ A_1(z)^{H^2} \cdot A_2(z)^{\chi(\mathcal O_S)} \cdot A_3(z)^{H\cdot K_S} \cdot A_4(z)^{K_S^2} \right ]=\binom{H\cdot K_S-n+1}{n}.\end{equation} 
The binomial expression \eqref{another} vanishes for the range 
\begin{equation}
\label{range}
n-1\leq k \leq 2n-1,
\end{equation}
covering in particular the vanishing \eqref{vv} in \cite{V}, and establishing Lehn's conjecture. The resulting closed formula on the $K3$ blowup,
$$\int_{S^{[n]}} s_{2n}(H^{[n]}) = \binom{H\cdot K_S-n+1}{n}\, \, \, \text{when} \, \, \, \chi (H) = 3n-1,$$
can also be seen geometrically: when maximally exploited, the Reider-type argument used by Voisin yields in fact the entire vanishing range \eqref{range} for the Segre class. This formula should be compared to the evaluation \eqref{clo} for $K3$ surfaces. However, unlike the $K3$ case where the Segre integrals were found for all values of $\chi$, the present closed expression holds conditionally on $\chi$ and $n$. 

Let us now establish \eqref{another}. Writing $H\cdot K_S=k$, we see that $$\chi(H)=3n-1 \implies H^2=k+6n-6\, .$$
We obtain $$a=H\cdot K_S+\chi(\mathcal O_S)=k+2, \,\,b=(H-K_S)^2+3\chi(\mathcal O_S)=-k+6n-1,$$ $$c=\chi(H)=3n-1.$$ Hence, in the light of \eqref{series}, we need to extract the coefficient of $z^n$ in the expression 
$$\frac{(1-w)^{k+2}(1-2w)^{-k+6n-1}}{(1-6w+6w^2)^{3n-1}}\, .$$ It is more convenient to express this coefficient as the residue $$\text{Res}_{z=0 } \,\omega$$ of the differential form $$\omega=\frac{(1-w)^{k+2}(1-2w)^{-k+6n-1}}{(1-6w+6w^2)^{3n-1}}\cdot\frac{dz}{z^{n+1}}.$$ Lehn's change of variables $$z=\frac{w(1-w)(1-2w)^{4}}{(1-6w+6w^2)^3}$$  is a nonsingular coordinate change near $w=0$ $$dz=\frac{(1-2w)^3}{(1-6w+6w^2)^3}dw\, .$$ Substituting, we obtain $$\omega=(1-w)^{k-n+1}(1-2w)^{-k+2n-2} \cdot \frac{dw}{w^{n+1}}\, .$$ A further change of variables $$w=\frac{u}{1+2u}$$ turns the form into $$\omega=(1+u)^{k-n+1}\cdot \frac{du}{u^{n+1}}\,.$$ The residue is now easily computed $$\text{Res}_{u=0} \,\omega=\binom{k-n+1}{n}\,,$$ thus confirming \eqref{another}. \qed
\vskip.1in
\noindent{\bf Remark.} Closed formulas for certain Segre integrals similar to \eqref{another} hold on blowups of all $K$-trivial surfaces. By the same methods it can shown that
\begin{itemize}
\item [(i)] If $S$ is the blowup of an Enriques surface at two points, then $$\int_{S^{[n]}}s_{2n}(H^{[n]})= \binom{H\cdot K_S-n+3}{n}$$ whenever $\chi(H)=3n-1$. 
\item [(ii)] If $S$ is the blowup of an abelian or bielliptic surface in three points, then $$\int_{S^{[n]}}s_{2n}(H^{[n]})= \binom{H\cdot K_S-n+5}{n}$$ whenever $\chi(H)=3n-1$. 
\end {itemize} 
\vspace{10pt}

\noindent{\bf Simpler form of the series.} For curves, the Segre series writes, according to \eqref{curve}, 
\begin{equation}\label{nwwn}
\sum_{n=0}^{\infty} z^n\int_{C^{[n]}} s_{n}({H}^{[n]})= A_1 (z)^d \cdot A_2 (z)^{\chi ({\mathcal O}_C)},
\end{equation} where
$$A_1(z)=1+t\, ,\,\,\,\,\,\, A_2(z)=\frac{(1+t)^2}{1+2t}\, $$
under the change of variables $$z=-t(1+t)\, .$$ 
For surfaces, a similar change of variables for surfaces simplifies the presentation of the universal series $A_1, A_2, A_3, A_4$ in \eqref{split} and is better suited for higher-rank generalizations. Specifically, 
setting \begin{equation}\label{changev}z=\frac{1}{2} t(1+t)^2, \text{ so that } w=\frac{1}{2}\left(1-\sqrt{\frac{1+t}{1+3t}}\right),\end{equation} a straightforward calculation using \eqref{series} 
yields:
\begin{eqnarray*}
A_1(z)&=&(1+t)^{\frac{1}{2}}\,,\\ 
A_2(z)&=&\frac{(1+t)^{\frac{3}{2}}}{(1+3t)^{\frac{1}{2}}}\,, \\
A_3(z)&=&\frac{\sqrt{1+t}+\sqrt{1+3t}}{2 (1+t)}\,, \\
 A_4(z)&=&\frac{4(1+t)^{\frac{1}{2}} (1+3t)^{\frac{1}{2}}} {\left (\sqrt{1+t}+\sqrt{1+3t}\right)^2}\, . 
 \end{eqnarray*}

\vspace{10pt}
\noindent{\bf Higher rank.} We discuss higher rank analogues of the above formulas. For a pair $(C, V)$ consisting of a nonsingular projective curve $C$ and a rank $r$ vector bundle $V$ of degree $d$, we have 
\begin{equation}\label{dxxd}
\sum_{n=0}^{\infty} z^n \int_{C^{[n]}} s_{n} (V^{[n]}) =A_1(z)^d \cdot A_2 (z)^{\chi(\mathcal O_C)},
\end{equation} for power series $A_1(z)$ and $A_2(z)$ depending  upon $r$. The series $A_1$ was conjectured in \cite {Wang}, though not in closed form, while the expression for $A_2$ was left open. Here, we prove the following result.
\begin {theorem} \label{hr} For formula \eqref{dxxd} in rank $r$, we have $$A_1(-t(1+t)^r)=1+t,\,\,\,   \ A_2(-t(1+t)^r)=\frac{(1+t)^{r+1}}{1+t(r+1)}\, .$$
\end{theorem}
\vskip.1in

\noindent{\it Proof of Theorem \ref{hr}.}
To find the series $A_1$ and $A_2$, we need only consider the projective line $C\simeq \mathbb P^1$ with the vector bundle 
$$V=\mathcal O_{\mathbb P^1}\otimes \mathbb C^{r-1} \oplus \mathcal O_{\mathbb P^1}(d)\, .$$ We obtain $$V^{[n]}=\mathcal O^{[n]}\otimes \mathbb C^{r-1}\oplus (\mathcal O(d))^{[n]}\, .$$ The Hilbert scheme of points is simply $(\mathbb P^1)^{[n]}\simeq \mathbb P^n$, and the universal subscheme $\mathcal Z\hookrightarrow  \mathbb P^n\times \mathbb P^1$ is given by $$\mathcal O(-\mathcal Z)=\mathcal O_{\mathbb P^n}(-1)\boxtimes \mathcal O_{\mathbb P^1}(-n)\, .$$ It follows that \begin{eqnarray*}\text{ch }\mathcal O(d)^{[n]}&=&\text{ch } {\mathbf R}\text{pr}_{\star}\left(\mathcal O_{\mathcal Z}\otimes \mathcal O_{\mathbb P^1}(d)\right)\\&=&\text{ch } {\mathbf R}\text{pr}_{\star}\left(\left(\mathcal O-\mathcal O(-\mathcal Z)\right)\otimes \mathcal O_{\mathbb P^1}(d)\right)\\&=&\text{ch}\left( H^0(\mathcal O_{\mathbb P^1}(d))\otimes \mathcal O_{\mathbb P^n}-H^{\bullet} (\mathcal O_{\mathbb P^1}(d-n))\otimes \mathcal O_{\mathbb P^n}(-1)\right)\\ &=& (d+1) - (d-n+1) \cdot \exp (-h)\end{eqnarray*} Here, we write $h$ for the hyperplane class on $\mathbb P^n$. We 
can then find the Chern roots of $(\mathcal O(d))^{[n]}$ yielding the following expression for the Segre class 
$$s(\mathcal O(d)^{[n]})=(1-h)^{d-n+1}\, .$$ Consequently $$s(V^{[n]})=(1-h)^{d-rn+r}\implies \int_{\mathbb P^n} s_n(V^{[n]})= (-1)^{n} \binom {d-rn+r}{n}\, .$$ 

We conclude that \begin{equation}\label{egl}\sum_{n=0}^{\infty} (-1)^n \binom{d-rn+r}{n} \cdot z^n= A_1(z)^d \cdot A_2(z).\end{equation} To finish the proof, we invoke the following result which was first proved in \cite {MOP} for $r=2$.

\setcounter{lemma}{2}
\begin{lemma} \label{l3}After the change of variables $$z=t(1+t)^r,$$ we have $$\sum_{n=0}^{\infty} \binom{d-rn+r}{n} \cdot z^n = \frac{(1+t)^{d+r+1}}{1+t(r+1)}.$$
\end{lemma} 

\proof We confirm that the coefficient of $z^n$ in the expression $$\frac{(1+t)^{d+r+1}}{1+t(r+1)}$$ equals $\binom{d-rn+r}{n}$ via a residue calculation. To this end, it suffices to prove that $$\text{Res}_{ z=0} \frac{(1+t)^{d+r+1}}{1+t(r+1)}\cdot \frac{dz}{z^{n+1}}=\binom{d-rn+r}{n}.$$ For the change of variables $z=t(1+t)^r$ we compute $$dz=(1+t)^{r-1}(1+t(r+1))\,dt.$$ Therefore, $$\text{Res}_{ z=0} \frac{(1+t)^{d+r+1}}{1+t(r+1)}\cdot \frac{dz}{z^{n+1}}=\text{Res}_{ t=0} \frac{(1+t)^{d-rn+r}}{t^{n+1}}\,dt=\binom{d-rn+r}{n}.$$\qed
 \vskip.1in

\noindent {\bf Surfaces.} For surfaces, a complete higher rank analogue of Lehn's conjecture is an open question. In this direction, several conjectures were recently formulated by D. Johnson \cite {J}, relating Segre theory to Verlinde theory in the Hilbert scheme context. Johnson's formulation of the conjectures was inspired by counts of points of $0$-dimensional Quot schemes and strange duality, similar to the strategy used to prove strange duality for curves in \cite {MO}. In the surface case, strange duality was pursued along these lines in \cite{bgj}. We sharpen the conjectures in \cite{J}, providing closed formulas for some of the series involved.

Specifically, consider a pair $(S, V)$ where $V$ is a rank $s$ vector bundle 
on a nonsingular projective surface $S$. The associated vector bundle $V^{[n]}$ on the Hilbert scheme has rank $sn$. By passing to resolutions, $V^{[n]}$ makes sense for all $K$-theory classes $V$. 

The following integrals of $V^{[n]}$ depend on five universal power series \begin{multline}\label{higherlehn}\sum_{n=0} z^n \int_{S^{[n]}} c_{2n}(V^{[n]})=\\ A_1(z)^{c_2(V)}\cdot A_2(z)^{\chi(c_1(V))}\cdot A_3(z)^{\frac{1}{2}\chi(\mathcal O_S)}\cdot A_4(z)^{K_S\cdot c_1(V)-\frac{1}{2} K_S^2} \cdot A_5(z)^{K_S^2}\, .\end{multline} After changing $V$ into $-V$ in $K$-theory, the above expressions turn into Segre integrals of higher rank vector bundles. Hence, equation \eqref{higherlehn} generalizes Lehn's formula. 

To connect with Verlinde theory, we recall first a result of \cite {EGL} regarding the holomorphic Euler characteristics of tautological line bundles: \begin{equation}\label{verlinde}\sum_{n=0}^\infty z^n \chi(S^{[n]}, H_{n}\otimes E^r)=f_r(z)^{\frac{1}{2}\chi(\mathcal O_S)} \cdot g_r(z)^{\chi(H)} \cdot a_r(z)^{H\cdot K_S-\frac{1}{2}K_S^2} \cdot b_r(z)^{K_S^2}.\end{equation} Here, $H_{n}$ denotes the line bundle induced by $H=\det V$ on the symmetric product, and $E$ is $-\frac{1}{2}$ of the exceptional divisor. By \cite {J} and \cite {EGL}, the two series corresponding to $K$-trivial surfaces are determined in closed form 
 $$f_r(z)=\frac{(1+t)^{r^2}}{1+r^2t},\,\,\,\,g_r(z)=1+t$$ after the change of variables $$z=t(1+t)^{r^2-1}.$$ As is usually the case, the series $a_r, b_r$ are unknown. 

Refining the conjectures{\footnote{The series $A_1,\ldots, A_5$ up to order 6 in $z$ were calculated in [J]. The numerical data in [J] 
played an important role in our formulation of Conjecture 1.}}
in \cite {J}, we provide closed expressions for the series in \eqref{higherlehn} corresponding to $K$-trivial surfaces. The last two series are surprisingly connected  in a {\it very precise} fashion to the unkown series $a_r, b_r$ of \eqref{verlinde}. 
\begin{conjecture} \label{pp23} Let $V\rightarrow S$ be a vector bundle{\footnote{For rank $s=1$  (corresponding to $r=0$), the Chern class
$c_{2n}(V^{[n]})$ 
is trivial for $n>0$ since $V^{[n]}$ is only of rank $n$.
The formulas of Conjecture \ref{pp23} are singular in the $r=0$ case.}} 
of rank $s=r+1$.
 After the change of variables $$z=-\frac{1}{r} t(1+t)^{-r}, \,\,\,\, w=\frac{t(-r+(-r+1)t)^{r^2-1}}{(-r(1+t))^{r^2}}\, ,$$ we have \begin{eqnarray*}A_1(z)&=&(-r)^{-r-1} \cdot (1+t)^{-r}\cdot (-r+(-r+1)t)^{r+1}\, , \\A_2(z)&=&(-r)^{r}\cdot (1+t)^{r-1} \cdot (-r+(-r+1)t)^{-r}\, , \\A_3(z)&=&(-r)^{r^2} \cdot (1+t-rt)^{-1}\cdot (1+t)^{(r-1)^2} \cdot (-r+t(-r+1))^{-r^2}\, ,\\A_4(z)&=&a_r(w)\, ,\\ A_5(z)&=&b_r(w)\, .\end{eqnarray*} 
\end {conjecture} 
Furthermore, using the solution of Lehn's conjecture, we are able to predict the first nontrivial{\footnote{We have $a_0=a_{\pm 1}=b_0=b_{\pm 1}=1\, .$}} 
examples
of the unknown series $a_r, b_r$ corresponding to $r=\pm 2$. 

\begin{conjecture} After the change of variables $$w=\frac{t(2+3t)^3}{16(1+t)^4}\, ,$$ we have \begin{eqnarray*}a_{-2}(w)&=&\frac{1}{a_2(w)}=\frac{2+3t}{\sqrt{1+t}} \cdot \frac{1}{\sqrt{1+t}+\sqrt{1+3t}}\, ,\\ b_{-2}(w)&=&b_2(w)=4 \sqrt{2+3t}\cdot \frac{(1+t)^{1/4}\cdot \sqrt{1+3t}}{(\sqrt{1+t}+\sqrt{1+3t})^{5/2}}\, .\end{eqnarray*} 
\end {conjecture} 
\noindent These formulas are connected to the series appearing in Lehn's rank $1$ formula.
 We have checked the term by term expansions pertaining to both $a_{\pm 2}$ and $b_{\pm 2}$ to high order. 

In case $S$ is a $K3$ surface, the series $a_r$ and $b_r$ 
play no role since $K_S$ vanishes. The expressions for the series $A_1, A_2, A_3$ were confirmed in \cite{mop2}, partially proving Conjecture \ref{pp23}.

\vspace{10pt}
\noindent {\bf Acknowledgements.} 
We thank C. Voisin and A. Szenes for discussions about Lehn's conjecture,
and
we thank D. Johnson for sharing with us additional numerical data.
A.M. was supported by 
the NSF through grant DMS 1601605. D. O. was supported by the NSF through grant  DMS 1150675. 
R.P. was supported by the Swiss National Science Foundation and
the European Research Council through
grants SNF-200020-162928 and ERC-2012-AdG-320368-MCSK.
R.P was also supported by SwissMap and the Einstein Stiftung in Berlin.

\end{document}